
\input amstex
\documentstyle{amsppt}
\magnification 1200
\loadbold
\TagsAsMath
\NoBlackBoxes

\define\PC{\Bbb P^1}
\define\Oh#1{\hat{#1}}
\define\pd#1#2{\dfrac{\partial#1}{\partial#2}}

\define\1{\sqrt{-1}\,}
\define\Ov#1{\overline{#1}}

\topmatter

\title On a conjecture of Varchenko
\endtitle

\author Roberto Silvotti
\endauthor

\affil Columbia University \endaffil

\date February 1995 \enddate

\dedicatory Dedicated to the memory of Aldo Andreotti.
\enddedicatory

\endtopmatter

\document

\head 1. Introduction \endhead

This note was motivated by a problem posed by Varchenko in \cite{7}. (See also
the earlier
paper by Aomoto \cite{2}). Let $f_I(x_1,\dots,x_n)$ ($I=1,\dots,N-1$) be linear
functions  of $n$
complex variables, and let $D_I=\{f_I=0\}\subset \Bbb C^n$ be the hyperplanes
defined as their
corresponding zero loci. The product
$\phi_{\lambda}=\prod_{I=1}^{N-1}\,f_I^{\lambda_I}$, where the $\lambda_I$ are
complex parameters,
is a multivalued holomorphic function on the complement $Y=\Bbb
C^n-\cup_{I=1}^{N-1}\,D_I$.

\proclaim{Varchenko's conjecture} If the exponents $\lambda_I$ are sufficiently
generic, then,
under certain broad conditions on the hyperplanes $D_I$:
\roster
\item The critical set of $\phi_{\lambda}$ in $Y$ is a union of isolated
points;
\item all critical points of $\phi_{\lambda}$ are nondegenerate and
\item their number is equal to $|\chi(Y)|$, the topological Euler
characteristic of $Y$ made
positive.
\endroster
\endproclaim

In its original version \cite{7}, the problem arose as a necessary step for
evaluating the
asymptotic behaviour of certain generalized hypergeometric integrals.
In the same paper, Varchenko also went on to prove the various assertions in
the case where the
family $\{D_1,\dots,D_{N-1}\}$ is a real arrangement, meaning that the linear
functions $f_I$ have
real coefficients. Subsequently, Orlik and Terao \cite{6} proved Varchenko's
conjecture in the
general situation of linear functions $f_I$ with complex coefficients. The
conditions are spelled out in \cite {6}: the family $\{D_1,\dots,D_{N-1}\}$
should be an
{\it essential arrangement\/}, i.e., the lowest dimensional multiple
intersections of hyperplanes
in the family should be isolated points.

The proofs in \cite{7}\cite{6} are of a combinatorial--topological nature and
thus
substantially rely on the assumption that the $D_I$ be hyperplanes in $\Bbb
C^n$. However, the
counting problem being purely topological, one would naturally expect that its
answer
should be found either in the evaluation of an appropriate characteristic class
or, alternatively, in a suitably construed Morse theoretical argument. In this
note we shall give
two independent proofs of a generalization of Varchenko's conjecture to the
case of
hypersurfaces in an algebraic manifold $X$. The first and most straightforward
proof is algebraic;
its main step consists in identifying and evaluating the number of critical
points as the top
Chern class of the sheaf of logarithmic 1--forms on a blow--up of $X$. The
second proof is an
application of  Morse theory; it is in part reminescent of the arguments used
by Andreotti and
Frankel \cite{1} and by Bott \cite{3} in their proofs of Lefschetz hyperplane
theorem. It also
bears some similarities to Aomoto's work \cite{2} concerning a naturally
related cohomology
of multivalued meromorphic forms on $X$. In either
case the assumption that the $D_I$ be hyperplanes becomes immaterial and can be
simply dispensed
with.

\medpagebreak

We shall thus consider the following generalization of Varchenko's problem. Let
$X$ be a smooth
projective variety of complex dimension $n$, that is a complex manifold which
can be embedded in
projective space $\Bbb P^m$, for some $m\ge n$, as the common zero set of
homogeneus polynomials.
Let also $D$ be a hypersurface in $X$ with (not necessarily smooth) irreducible
components
$D_1$,\dots,$D_N$. We shall consider nowhere vanishing multivalued
holomorphic sections of a flat complex line bundle on $Y=X-D$
required to have ``power behaviour'' near $D$. Concretely, such a section
is a multivalued, holomorphic and nowhere vanishing function $\phi_{\lambda}$
on $Y$ with the
following property: If $\{f_I=0\}$ is a local defining equation of $D_I$ on a
sufficiently small
neighborhood $U\subset X$ of a smooth point $p\in D$, $\phi_{\lambda}$ has the
local form
$$
\phi_{\lambda}|_U=f_I^{\lambda_I}\,h,
$$
where $\lambda_I$ is a complex number and $h$ is holomorphic  and
non--vanishing throughout
$U$. It can be readily verifed that the number $\lambda_I$ neither depends on
the
chosen smooth point $p$ nor on the choice of the local defining function $f_I$
near $p$; we shall
call it the {\it order\/} of $\phi_{\lambda}$ along $D_I$. The space of orders
$\lambda=(\lambda_1,\dots,\lambda_N)$ is generally a homogeneus hyperplane
$\Lambda$ in
$C^N$ defined by a linear polynomial with positive integer coefficients.

The following preliminary proposition is proven in Section 2.

\proclaim{Proposition 1.1} Let $\phi_{\lambda}$ be as above, and assume that
there is a dense open
subset $V\subset \Lambda$ such that, for $\lambda\in V$, the critical points of
$\phi_{\lambda}$
in $Y$ are all non--degenerate. Then there is a linear homogeneus subset
$A\subset\Lambda$ such that
the number of critical points is independent of $\lambda\in V-A\cap V$.
\endproclaim

Our main theorem then solves the counting problem.

\proclaim{Theorem 1.2} Under the same assumption, then, for $\lambda\in V-A\cap
V$, the number
of critical points of $\phi_\lambda$ in $Y$ is equal to $(-1)^n\,\chi(Y)$, the
topological Euler
characteristic of $Y$ up to a sign.
\endproclaim

\remark{Remark 1.3} The subset $A$ has the following
geometric origin. By Hironaka's resolution of singularities theorem, there
exists a resolution
$\sigma\:\hat X\to X$ in which the preimage $\hat Y=\hat X-\hat D$ of $Y$
becomes the complement of
a normal crossing divisor $\hat D$ and such that $\sigma$ restricts to an
isomorphism from $\hat Y$
to $Y$. Let $\{g_i=0\}$ be a local defining equation of the component $\hat
D_i$ of $\hat D$. Near a
smooth point of $\hat D_i$, the pull--back $\sigma^*\phi_{\lambda}$ has the
local form
$g_i^{\hat\lambda_i}\,h$, where $h$ is holomorphic and nowhere vanishing and
the
order $\hat\lambda_i\equiv\hat\lambda_i(\lambda)$ of $\sigma^*\phi_{\lambda}$
along $\hat D_i$ is
a homogeneus linear polynomial in the original orders $\lambda_I$ with positive
integer
coefficients. The subset $A$ is then defined as the the zero set
$$
A =\left\{\lambda\in\Lambda\mid\prod_{i\in \hat
I}\,\hat\lambda_i(\lambda)=0\right\},
$$
where $\hat I$ indexes the irreducible components of $\hat D$.

In fact, one can easily obtain (see Section 4) a general---though less
explicit---formula for the
number of critical points valid for \underbar{any} $\lambda\in V$.
\endremark

\bigpagebreak

\remark{Remark 1.4} When one allows for possibly degenerate but still isolated
critical points,
the result still holds true with the following modification (see Section 4).
For $U_p$ a small
neighborhood of a critical point $p$ of $\phi_\lambda$, let
$\varphi_\lambda|_{U_p^*}\:U_p^*=U_p-\{p\}\to\Bbb C^n-\{0\}$ denote the map
whose components are
the components of the 1--form $d\log\phi_\lambda$ on $U_p$. Then one has
$$
\chi(Y)=(-1)^n\,\sum_{\{\text{critical points $p$}\}}\,\,
\operatorname{deg}\varphi_\lambda|_{U_p^*},
$$
where $\operatorname{deg}\varphi_\lambda|_{U_p^*}$ denotes the topological
degree.
\endremark

\bigpagebreak

In this note we shall not address the question of finding ``effective''
geometric conditions on the
hypersurfaces $D_I$ in order for an open set $V\subset\Lambda$ satisfying the
above hypotheses to
exist. One can easily find, however, large classes of non--trivial examples.
Here we list two.

\subhead Example 1.5\endsubhead  This is the case of Varchenko--Orlik--Terao.
Let $X=\Bbb P^n$
and let $D_N$ be the hyperplane at infinity. The remaining $D_I$ are the
hyperplanes
$\{f_I=0\}\subset\Bbb C^n$, where the linear functions
$f_I=\sum_{i=1}^n\,a_{iI}x_i+b_I$ are such
that the $n\times (N-1)$ constant matrix $(a_{iI})$ has rank equal to $n$. The
family
$D_1,\dots,D_{N-1}$ is thus an affine essential arrangement. In this case Orlik
and Terao
\cite{6: Section 4} have shown that the conditions of Proposition 1.1 and
Theorem 1.2 are satisfied.

\subhead Example 1.6\endsubhead This example should be contrasted with Example
1.5. It is meant to
give a simple illustration of how the geometric conditions on the intersections
that have to be
imposed on the $D_I$ in the case of hyperplanes (i.e., that the arrangement be
essential), become
superfluous when at least one of the $D_I$ is a hypersurface of higher degree.

Let $D=D_1\cup D_2$ be the hypersurface in $\Bbb P^n$ whose components are
respectively the zero
locus of the degree $d\ge 2$ homogeneus polynomial $F_1=\sum_{i=0}^n\,X_i^d$
and of $F_2=X_0$.
Thus, $D_2$ being the hyperplane at infinity, $Y=\Bbb P^n-D$ is the complement
in $\Bbb C^n=\Bbb
P^n-D_2$ of the affine hypersurface $D_1\cap\Bbb
C^n=\{f_1=\sum_{i=1}^n\,x_i^d+1=0\}$,
where $x_i=X_i/X_0$ are affine coordinates on $\Bbb C^n$. The multivalued
function
$\phi_\lambda=f_1^\lambda$ depends on a single parameter $\lambda\in\Bbb C$.
The critical set
$C_\lambda$ of $\phi_\lambda$, defined by the equations $\lambda
x_1=\dots=\lambda x_n=0$, consists
of a single point, the origin $C_\lambda=\{0\}$ in $\Bbb C^n$, for
$\lambda\in\Bbb C^*$, and of all
of $Y$, $C_0=Y$, for the special value $\lambda=0$. From the Hessian matrix
$$
\operatorname{Hess}(\phi_\lambda) =\lambda\,d(d-1)\,
\operatorname{diag}\bigl(x_1^{d-2},\dots,x_n^{d-2}\bigr),
$$
one sees that the critical point is degenerate unless $d=2$. On the other hand,
let $U_0$ be a
small neighborhood of the origin in $\Bbb C^n$; the topological degree of the
map
$\varphi_\lambda|_{U_0^*}\:U_0^*=U_0-\{0\}\to\Bbb C^n-\{0\}$ sending $x$ to
$\bigl(\partial_{x_1}\log\phi_\lambda(x),\dots,\partial_{x_1}\log\phi_\lambda(x)\bigr)=
\lambda d \,f_1^{-1}\,(x_1^{d-1},\dots,x_n^{d-1})$, is equal to $(d-1)^n$.
Theorem 1.2 and the
following remark say in this case that $\chi(Y)=\chi(\Bbb C^n)-\chi(D_1\cap\Bbb
C^n)=(-1)^n\,(d-1)^n$. This is in agreement with the well--known fact that the
cohomology
of the affine hypersurface $D_1\cap\Bbb C^n$ is non--vanishing only in degree
$n-1$ and $2n-2$,
where $\operatorname{dim}H^{n-1}(D_1\cap\Bbb C^n)=(d-1)^n$,
$\operatorname{dim}H^{2n-2}(D_1\cap\Bbb C^n)=1$.

\remark{Acknowledgements} The author was partially supported by NSF grant DMS
92--04196.
He has benefitted from discussions with R. Friedman,
Y. Karpishpan, J. Morgan, D.H. Phong, H. Pinkham and S. Wu. He would
like to especially thank P. Orlik and H. Terao for several exchanges, and
Nicholas
Shepherd--Barron for extremely helpful indications.
\endremark

\head 2. A few preliminary observations \endhead

The purpose of the following observations is to show that the number
of critical points of $\phi_\lambda$ can be easily identified, for almost all
$\lambda$, with a
topological invariant of the pair $(X,D)$.

It will be more convenient to work with single--valued objects on
$Y$ rather than directly with $\phi_\lambda$. We thus note first of all that,
since
$\phi_\lambda$ is nowhere vanishing on $Y$, the critical set $C_\lambda$ of
$\phi_\lambda$ is
precisely the set of points in $Y$ where the meromorphic 1--form
$d\log\phi_\lambda=\frac
1{\phi_\lambda}d\phi_\lambda$ vanishes. Let
$\varphi_{\lambda,i}=\dfrac{\partial}{\partial x_i}\log\phi_\lambda\:U\to\Bbb
C$ for
$i=1,\dots,n$ denote the components of $d\log\phi_\lambda$ on any open set $U
\subset Y$ of a
coordinate cover of $Y$ with local coordinate $x=(x_1,\dots,x_n)\:U\to\Bbb
C^n$. Then
$C_\lambda$ is the analytic subvariety of $Y$ with local defining equations
$\varphi_{\lambda,1}=\dots=\varphi_{\lambda,n}=0$ on $U$. Second, note that the
Hessian
matrix of $\phi_\lambda$ on $U$ is related to the Jacobian of the map
$\varphi_\lambda\:U\to
\Bbb C^n$, $p\mapsto
\bigl(\varphi_{\lambda,1}(p),\dots,\varphi_{\lambda,n}(p)\bigr)$:
$$
\align
\dfrac{\partial^2\phi_\lambda}{\partial x_j\partial x_i} &=
\phi_\lambda\,\dfrac{\partial \varphi_{\lambda,i}}{\partial x_j}+
\phi_\lambda\,\varphi_{\lambda,i}\,\varphi_{\lambda,j}\\
&=\phi_\lambda\,\operatorname{Jac}(\varphi_\lambda)_{ij}+
\phi_\lambda\,\varphi_{\lambda,i}\,\varphi_{\lambda,j}.
\endalign
$$
Again, since $\phi_\lambda$ is never zero on $Y$, it follows that a critical
point $p\in
C_\lambda$ is non--degenerate if and only if the
determinant $\operatorname{det}\operatorname{Jac}(\varphi_\lambda)(p)\ne 0$. As
usual, we say in
this case that $p$ is a {\it non--degenerate zero\/} of $d\log\phi_{\lambda}$.

\medpagebreak

We now turn to the proof of of Proposition 1.1. If $\hat X @>\sigma>> X$ is the
blow--up of Remark
1.2, let $\hat D=\sigma^{-1}(D)$ and $\hat Y=\sigma^{-1}(Y)$. Since the
Jacobian of $\sigma$ is a
holomorphic non--singular matrix on $\hat Y$, one may easily verify that the
1--form
$d\log\sigma^*\phi_{\lambda}\in H^0\bigl(\hat X,\Omega_{\hat X}^1(*\hat
D)\bigr)$ has a
non--degenerate zero at $\hat p\in\hat Y$ if and only if
$\hat p=\sigma^{-1}(p)$, where $p\in Y$ is a non--degenerate zero of
$d\log\phi_{\lambda}$.
Thus clearly, if $\lambda\in V$, the cardinality of $C_\lambda$ is given by
$$
\bigl(\text{\# of critical points of $\phi_{\lambda}$ on $Y$}\bigr) =
\bigl(\text{\# of zeroes of $d\log\sigma^*\phi_{\lambda}$ on $\hat Y$}\bigr).
$$
Let us recall that the sheaf of logarithmic 1--forms $\Omega_{\hat X}^1(\log
\hat D)$
is the the sheaf of those meromorphic 1--forms $\omega$ on $\hat X$ which are
holomorphic on $\hat
X-\hat D$ and have the following local property near $\hat D$: For any small
open neighborhood
$U\subset\hat X$ of $\hat D$ on which $\hat D$ has a local defining equation
$g=0$, both
$g\,\omega$ and $g\,d\omega$ are holomorphic throughout $U$. Also, recall the
algebraic subset
$A\subset\Lambda$ introduced in Remark 1.2.

\proclaim{Lemma 2.1} Let $\lambda\in\Lambda$. Then: (i)\ The 1--form
$d\log\sigma^*\phi_{\lambda}$
is an element of $\Gamma\bigl(\hat X,\Omega_{\hat X}^1(\log \hat D)\bigr)$.
(ii)\ For $\lambda\in
\Lambda - A$, $d\log\sigma^*\phi_{\lambda}$ has a pole along \underbar{every}
component of $\hat D$.
\endproclaim

\demo{Proof} The assertions being local, it suffices to consider the
pull--back $\sigma^*\phi_{\lambda}$ of $\phi_{\lambda}$ on an arbitrarily small
neighborhood
$U\subset \hat X$ of a point where exactly $m$ components of $\hat D$, say
$\hat D_1,\dots,\hat
D_m$, intersect. If $\{g_i=0\}$ is a local defining equation of $\hat D_i$ on
$U$,
$\sigma^*\phi_{\lambda}$ has the local form $g_1^{\hat\lambda_1}\cdots
g_m^{\hat\lambda_m}\,h$,
where $h$ is holomorphic and nowhere vanishing throughout $U$. Assertion {\sl
(i)\/} is thus
self--evident. Moreover, if none of the $\hat\lambda_i$ is zero---i.e., if
$\lambda\notin A$---then
$d\log\sigma^*\phi_{\lambda}|_U$ has a logarithmic pole along $\hat D_i$ for
all $i=1,\dots,m$,
and {\sl (ii)\/} is also clear.
\enddemo

If $\lambda\in V-A\cap V$, the number of zeroes of
$d\log\sigma^*\phi_{\lambda}$ is
therefore a topological invariant of the vector bundle $\Omega_{\hat X}^1(\log
\hat D)$ on $\hat
X$. In particular, it is obviously independent of $\lambda$.

\head 3. Gauss--Bonnet for the complement of a divisor\\
         (First proof of Theorem 1.2)                    \endhead

For $\lambda\in V-A\cap V$, the number of zeroes of
$d\log\sigma^*\phi_{\lambda}$---and hence the
number of critical points of $\phi_\lambda$ in $Y$---has just been identified
with the top Chern
number, $\int_X c_n\bigl(\Omega_{\hat X}^1(\log \hat D)\bigr)$, of
$\Omega_{\hat X}^1(\log \hat D)$.
(For this standard interpretation of the top Chern class of a holomorphic
vector bundle see e.g.
\cite{4: section 3 of Chapter 3}). By Theorem 4.1 below,
this coincides with the topological Euler
characteristic of $\hat Y$ up to a sign factor of $(-1)^n$. Since $\hat Y$ is
isomorphic to $Y$,
then $\chi(\hat Y)=\chi(Y)$, which  concludes our first proof of Theorem 1.2.

\proclaim{Theorem 4.1 \footnote{The author is indebted to N. Shepherd--Barron
for having pointed out
this version of the Gauss--Bonnet formula to him. However, the author has  been
unable to locate a
proof anywhere in the literature.}} Let $D=\sum_{I=1}^N\,D_I$ be a normal
crossing divisor in a
smooth projective variety $X$ of complex dimension $n$ and let $\Omega_X^1(\log
D)$ be the rank $n$
holomorphic vector bundle on $X$ whose sections are the 1--forms on $X$ with
logarithmic poles along
$D$. Then the top Chern number of $\Omega_X^1(\log D)$ is given by the Euler
characteristic of the
complement, $\chi(X-D)=\sum_{i=1}^n\,\operatorname{dim}_{\Bbb C}H^i(X-D,\Bbb
C)$, up to a sign,
$$
\int_X c_n\bigl(\Omega_X^1(\log D)\bigr)=(-1)^n\,\chi(X-D).
$$
\endproclaim

\demo{Proof} We shall first express the Chern classes of $\Omega_X^1(\log D)$
in terms of the Chern
classes of the holomorphic cotangent bundle $\Omega_X^1=T_X^*$ and of the line
bundles $[D_I]$
associated with the various components $D_I$. For $E$ a holomorphic vector
bundle on $X$, we shall
denote its total Chern class by $c(E)=\sum_{i\ge 0}\,c_i(E)$ with the usual
convention
$c_0(E)\equiv 1$. The Poincar\'e residue map gives the exact sequence of
sheaves on $X$
$$
0@>>> \Omega_X^1@>>>\Omega_X^1(\log D)@>\operatorname{residue}>>\Cal O_{\tilde
D}=\bigoplus_{I=1}^N\,
\Cal O_{D_I}@>>>0,
$$
where the $\Cal O_{D_I}$ are to be viewed as the sheaves on $X$ extending the
structure sheaves
$\Cal O_{D_I}$ by zero outside the divisors $D_I$. The resulting identity among
total Chern classes,
$c\bigl(\Omega_X^1(\log D)\bigr)=c(\Omega_X^1)\,c(\Cal O_{\tilde D})=
c(\Omega_X^1)\,\prod_{I=1}^N\,c(\Cal O_{D_I})$, gives
$$
\align
c_n\bigl(\Omega_X^1(\log D)\bigr)=&
\sum_{i=0}^n\,\,\sum_{j_1+\dots+j_N=i}\,c_{n-i}(\Omega_X^1)\,c_{j_1}(\Cal
O_{D_1})\cdots
c_{j_N}(\Cal O_{D_N}) \\
&=c_n(\Omega_X^1)+
\sum_{i=1}^n\,\,\sum_{j_1+\dots+j_N=i}\,c_{n-i}(\Omega_X^1)\,c_{j_1}(\Cal
O_{D_1})\cdots c_{j_N}(\Cal O_{D_N}).
\endalign
$$
Moreover, one has for every $I$ an exact sequence of sheaves on $X$,
$$
0@>>>\Cal O_X\bigl([-D_I]\bigr)=[-D_I]@>>>\Cal
O_X@>\operatorname{restriction}>>\Cal O_{D_I}@>>>0,
$$
implying $c(\Cal O_X)=1=c(\Cal O_{D_I})\,c\bigl([-D_I]\bigr)=
c(\Cal O_{D_I})\,\left(1+c_1\bigl([-D_I]\bigr)\right)$. Thus
$$
c(\Cal O_{D_I})=\sum_{j\ge 0}\,(-1)^j c_1\bigl([-D_I]\bigr)^j=
\sum_{j\ge 0}\, c_1\bigl([D_I]\bigr)^j,
$$
and the Chern classes of $\Cal O_{D_I}$ are $c_j(\Cal
O_{D_I})=c_1\bigl([D_I]\bigr)^j$.
Now, on the one hand we have the Gauss--Bonnet formula $\int_X
c_n(T_X^*)=(-1)^n
\int_X c_n(T_X)=(-1)^n\chi(X)$ for the top Chern class of $\Omega_X^1$; on the
other hand, by a
standard excision argument, we have the addition formula
$\chi(X-D)=\chi(X)-\chi(D)$. Hence
$$
\int_X c_n\bigl(\Omega_X^1(\log D)\bigr)=(-1)^n\,\chi(X) +
\sum_{i=1}^n\,\,\sum_{j_1+\dots+j_N=i}\,\int_X c_{n-i}(\Omega_X^1)\,
                     c_1\bigl([D_1]\bigr)^{j_1}\cdots
c_1\bigl([D_N]\bigr)^{j_N}
$$
and the sought for result is equivalent to the

\proclaim{Claim} Let $D_1,\dots, D_N$ be smooth and normally intersecting
divisors in a
$n$--dimensional smooth projective variety $X$, and let
$D=\bigcup_{I=1}^N\,D_I$. Then the following
identity holds
$$
(-1)^{n-1}\chi(D) =\sum_{i=1}^n\,\,\sum_{j_1+\dots+j_N=i}\,\int_X
c_{n-i}(\Omega_X^1)\,
                     c_1\bigl([D_1]\bigr)^{j_1}\cdots
c_1\bigl([D_N]\bigr)^{j_N}.
\tag 1
$$
\endproclaim

In order to prove the claim, let us first consider a smooth divisor $D_I$ in
$X$. Since $D_I$ has
complex codimension 1 in $X$, the normal bundle $N_{D_I/X}$ is a line bundle on
$D_I$ equivalent to
the restriction $[D_I]|_{D_I}$. The $C^\infty$ decomposition
$\left.T_X\right|_{D_I}=T_{D_I}\oplus
N_{D_I/X}=T_{D_I}\oplus [D_I]|_{D_I}$ implies the identity
$c(T_X)|_{D_I}=c\bigl(\left.T_X\right|_{D_I}\bigr)=c(T_{D_I})\,c\bigl([D_I]|_{D_I}\bigr)=c(T_{D_I})\,
c([D_I])|_{D_I}$ of Chern polynomials. The resulting equalities of cohomology
classes on $D_I$,
$$
c_i(T_X)|_{D_I}=c_i(T_{D_I})+c_{i-1}(T_{D_I})\,c_1([D_I])|_{D_I}
\qquad\text{for $i=1,\dots,n-1$}
$$
together with the general relations $c_i(E)=(-1)^i c_i(E^*)$ between the Chern
classes of a vector
bundle and of its dual, give
$$
c_i(T_X^*)|_{D_I}=c_i(\Omega_X^1)|_{D_I}=c_i(T_{D_I}^*)-c_{i-1}(T_{D_I}^*)\,c_1([D_I])|_{D_I}
\qquad\text{for $i=1,\dots,n-1$}.
\tag 2
$$

Our proof of the claim now proceeds by induction on the number $N$ of divisors.

\smallpagebreak

\par\noindent {\it Step 1.\/}\ For $N=1$, let $D=D_1$. The claim follows at
once from \thetag 2:
$$
\align
\sum_{i=1}^n\,\int_X c_{n-i}(\Omega_X^1)\,c_1([D])^i &=
\sum_{i=1}^n\,\int_D c_{n-i}(\Omega_X^1)\,c_1([D])^{i-1} \\
&= \sum_{i=1}^n\,\int_D c_{n-i}(T_D^*)\,c_1([D])^{i-1}-
 \sum_{i=1}^{n-1}\,\int_D c_{n-1-i}(T_D^*)\,c_1([D])^i \\
&=\sum_{i=0}^{n-1}\,\int_D c_{n-1-i}(T_D^*)\,c_1([D])^i-
 \sum_{i=1}^{n-1}\,\int_D c_{n-1-i}(T_D^*)\,c_1([D])^i\\
&=\int_D c_{n-1}(T_D^*)=(-1)^{n-1}\chi(D),
\endalign
$$
where in the first step we have used the fact that $c_1([D])$ is Poincar\'e
dual to the
fundamental class of $D$.

\smallpagebreak

\par\noindent {\it Step 2.\/}\ For general $N>1$ we decompose the sum in
\thetag 1 into the terms
with $j_N=0$ and those with $j_N\ge 1$:
$$
\multline
\sum_{i=1}^n\,\,\sum_{j_1+\dots+j_N=i}\,c_{n-i}(\Omega_X^1)\,
                     c_1\bigl([D_1]\bigr)^{j_1}\cdots
c_1\bigl([D_N]\bigr)^{j_N}=\\
\sum_{i=1}^n\,\,\sum_{j_1+\dots+j_{N-1}=i}\,c_{n-i}(\Omega_X^1)\,
                     c_1\bigl([D_1]\bigr)^{j_1}\cdots
c_1\bigl([D_{N-1}]\bigr)^{j_{N-1}}+\\
\sum_{i=1}^n\,\,\sum\Sb j_1+\dots+j_N=i\\ j_N\ge 1\endSb
\,c_{n-i}(\Omega_X^1)\,
                     c_1\bigl([D_1]\bigr)^{j_1}\cdots
c_1\bigl([D_N]\bigr)^{j_N}.
\endmultline
$$
The second sum on the right hand side can be computed using \thetag 2,
$c_i(\Omega_X^1)|_{D_N}=c_i(\Omega_{D_N}^1)-c_{i-1}(\Omega_{D_N}^1)\,c_1([D_N])|_{D_N}$:
$$
\multline
\sum_{i=1}^n\,\,\sum\Sb j_1+\dots+j_N=i\\ j_N\ge 1\endSb \,\int_X
c_{n-i}(\Omega_X^1)\,c_1\bigl([D_1]\bigr)^{j_1}\cdots
c_1\bigl([D_N]\bigr)^{j_N}=\\
\shoveleft{\sum_{i=1}^n\,\,\sum\Sb j_1+\dots+j_N=i\\ j_N\ge 1\endSb
\,\int_{D_N}
c_{n-i}(\Omega_X^1)\,c_1\bigl([D_1]\bigr)^{j_1}\cdots
c_1\bigl([D_N]\bigr)^{j_N-1}=}\\
\shoveleft{\sum_{i=1}^n\,\,\sum\Sb j_1+\dots+j_N=i\\ j_N\ge 1\endSb
\,\int_{D_N}
c_{n-i}(\Omega_{D_N}^1)\, c_1\bigl([D_1]\bigr)^{j_1}\cdots
c_1\bigl([D_N]\bigr)^{j_N-1}-}\\
\shoveleft{\,\sum_{i=1}^{n-1}\,\,\sum\Sb j_1+\dots+j_N=i\\ j_N\ge 1\endSb
\,\int_{D_N}
c_{n-1-i}(\Omega_{D_N}^1)\,c_1\bigl([D_1]\bigr)^{j_1}\cdots
c_1\bigl([D_N]\bigr)^{j_N}=}\\
\shoveleft{\int_{D_N}c_{n-1}(\Omega_{D_N}^1)+
\sum_{i=1}^{n-1}\,\,\sum_{j_1+\dots+j_{N-1}=i} \,\int_{D_N}
c_{n-1-i}(\Omega_{D_N}^1)\,c_1\bigl([D_1]\bigr)^{j_1}\cdots
c_1\bigl([D_{N-1}]\bigr)^{j_{N-1}}}.
\endmultline
$$
We have thus found the recursion
$$
\multline
\sum_{i=1}^n\,\,\sum_{j_1+\dots+j_N=i}\,\int_X c_{n-i}(\Omega_X^1)\,
                     c_1\bigl([D_1]\bigr)^{j_1}\cdots
c_1\bigl([D_N]\bigr)^{j_N}=\\
\sum_{i=1}^n\,\,\sum_{j_1+\dots+j_{N-1}=i}\,\int_X c_{n-i}(\Omega_X^1)\,
                     c_1\bigl([D_1]\bigr)^{j_1}\cdots
c_1\bigl([D_{N-1}]\bigr)^{j_{N-1}}+
                     (-1)^{n-1}\,\chi(D_N)+\\
\sum_{i=1}^{n-1}\,\,\sum_{j_1+\dots+j_{N-1}=i} \,\int_{D_N}
c_{n-1-i}(\Omega_{D_N}^1)\,c_1\bigl([D_1]\bigr)^{j_1}\cdots
c_1\bigl([D_{N-1}]\bigr)^{j_{N-1}}.
\endmultline
$$
Let us now assume that the claim is true for a number $N-1$ of smooth divisors
with transverse
intersections. This means that the first term above is given by
$$
\sum_{i=1}^n\,\,\sum_{j_1+\dots+j_{N-1}=i}\,\int_X c_{n-i}(\Omega_X^1)\,
                     c_1\bigl([D_1]\bigr)^{j_1}\cdots
c_1\bigl([D_{N-1}]\bigr)^{j_{N-1}}=
(-1)^n\,\chi\left(\cup_{I=1}^{N-1}\,D_I\right);
$$
moreover, since the $N-1$ divisors $D_1\cap D_N$,\dots, $D_{N-1}\cap D_N$ in
the compact
$(n-1)$--dimensional smooth projective variety $D_N$ are smooth and have
transverse intersections,
by the induction hypothesis we also have
$$
\multline
\sum_{i=1}^{n-1}\,\,\sum_{j_1+\dots+j_{N-1}=i} \,\int_{D_N}
c_{n-1-i}(\Omega_{D_N}^1)\,c_1\bigl([D_1]\bigr)^{j_1}\cdots
c_1\bigl([D_{N-1}]\bigr)^{j_{N-1}}=\\
(-1)^{n-2}\,\chi\left(\cup_{I=1}^{N-1}\,D_I\cap D_N\right).
\endmultline
$$
But the Mayer--Vietoris cohomology exact sequence of the pair
$\left(\cup_{I=1}^{N-1}\,D_I,\,D_N\right)$ implies the relation
$$
\chi(D)=\chi\left(\cup_{I=1}^N\,D_I\right)=\chi\left(\cup_{I=1}^{N-1}\,D_I\right)+\chi(D_N)-
\chi\left(\cup_{I=1}^{N-1}\,D_I\cap D_N\right)
$$
among Euler characteristics. It follows that
$$
\multline
\sum_{i=1}^n\,\,\sum_{j_1+\dots+j_N=i}\,\int_X c_{n-i}(\Omega_X^1)\,
                     c_1\bigl([D_1]\bigr)^{j_1}\cdots
c_1\bigl([D_N]\bigr)^{j_N}=\\
(-1)^{n-1}\,\chi\left(\cup_{I=1}^{N-1}\,D_I\right)+(-1)^{n-1}\,\chi(D_N)+
(-1)^{n-2}\,\chi\left(\cup_{I=1}^{N-1}\,D_I\cap D_N\right)= \\
(-1)^{n-1}\,\chi(D),
\endmultline
$$
and the proof of the induction step is complete.

\enddemo

\head 4. Variants  \endhead

In this brief section we discuss the variants mentioned in Remark 1.3 and
Remark 1.4.

\subhead 4.1\endsubhead  Even if $\lambda\in A\cap V$ one can still in
principle compute the number
of critical points as follows.  For $\lambda\in\Lambda$, let us thus introduce
the divisor $\hat
D(\lambda)$ given by those irreducible components of $\hat D$ along which the
order of
$\sigma^*\phi_{\lambda}$ is non-zero,
$$
\hat D(\lambda)=\bigcup_{\hat\lambda_i\ne 0}\,\hat D_i \,\, \subset \hat D.
$$
Note that $\hat D(\lambda)=\hat D$ if and only if $\lambda\in \Lambda-A$. One
can immediately
deduce the following sharper version of Lemma 2.1:

\proclaim\nofrills{} Let $\lambda\in\Lambda$. Then the 1--form
$d\log\sigma^*\phi_{\lambda}$ is
an element of $\Gamma\bigl(\hat X,\Omega_{\hat X}^1(\log \hat
D(\lambda))\bigr)$ having a pole along
\underbar{every} component of $\hat D(\lambda)$.
\endproclaim

By the same argument as above, one thus deduces the more general formula
$$
\bigl(\text{\# of critical points of $\phi_{\lambda}$ on $Y$}\bigr)=
(-1)^n\,\chi\bigl(\hat X-\hat D(\lambda)\bigr)
$$
for any $\lambda\in V$. The practical usefulness of this formula for producing
numerical
predictions may however be limited to those concrete situations where one can
easily relate the
Euler characteristic of $\hat X-\hat D(\lambda)$ to the topology of the
blow--down of $\hat
D(\lambda)$.

\subhead 4.2\endsubhead Let $\lambda\in V$. We shall allow the section
$d\log\sigma^*\phi_{\lambda}\in \Gamma\bigl(\hat X,\Omega_{\hat X}^1(\log \hat
D(\lambda))\bigr)$
to have isolated but possibly degenerate zeroes $p$. The top Chern class
$c_n\bigl(\Omega_{\hat X}^1(\log \hat D(\lambda)\bigr)$ is Poincar\'e dual to
the degeneracy
cycle $\sum_{p}\,m_p\,p$. Here the multiplicity $m_p$ is the intersection
number at $p$ of
the $n$ divisors in $\hat Y$ having the local defining equations
$\hat\varphi_{\lambda,i}=\partial_{x_i} \log\sigma^*\phi_{\lambda}=0$.
Equivalently, $m_p$ is
the topological degree of the map
$\hat\varphi_{\lambda}|_{U_p^*}\:U_p^*=U_p-\{p\}\to\Bbb
C^n-\{0\}$ with components $\hat\varphi_{\lambda,i}|_{U_p^*}$, $U_p$ being a
small neighborhood of
$p$. It follows that
$$
(-1)^n\,\chi\bigl(\hat X-\hat D(\lambda)\bigr)=\int_{\hat
X}c_n\bigl(\Omega_{\hat X}^1(\log \hat
D(\lambda)\bigr)=\sum_p\,m_p.
$$
If $\lambda\in V-A\cap V$, this is the formula given in Remark 1.4.

\head 5. Morse theoretic proof of Theorem 1.2 \endhead

The main idea of this second proof consists---loosely speaking---of
interpreting
$|\phi_\lambda|^2=\phi_\lambda\,\phi_\lambda^*$ as a Morse function defined on
the
submanifold of $X$ obtained from this latter by deleting a tiny neighborhood of
the
hypersurface $D$. More precisely, we shall actually work in the blow--up $\hat
X @>\sigma>> X$ in
which $\hat Y@>\sim>> Y$ is realized as the complement of a normal crossing
divisor $\hat D$.
For simplicity, we shall denote the pull--back $\sigma^*\phi_\lambda$ by
$\hat\phi_\lambda$.
As already observed, $\hat\phi_\lambda$ has a non--degenerate critical point at
$p\in\hat Y$ if
and only if $p$ is the preimage of a non--degenerate critical point of
$\phi_\lambda$ in $Y$.
Let $x_1,\dots,x_n$ and $v_1,\dots,v_n,w_1,\dots, w_n$ be,
respectively, analytic and real local coordinates on $\hat X$ with $x_i=v_i+\1
w_i$. One
immediately verifies that the critical set of the real--valued function
$|\hat\phi_\lambda|^2$,
i.e., the set of points in $\hat Y$ where all partial derivatives
$\partial_{v_i}\,|\hat\phi_\lambda|^2$, $\partial_{w_i}\,|\hat\phi_\lambda|^2$
vanish, coincides
with the critical set of $\hat \phi_\lambda$ in $\hat Y$. The Hessian of $|\hat
\phi_\lambda|^2$
at a critical point $p$ is given by the $2n\times 2n$ matrix
$$
\align
\operatorname{Hess}(|\hat\phi_\lambda|^2)(p) & =
\pmatrix
\left(\partial_{v_i}\partial_{v_j}\,|\hat\phi_\lambda|^2\right)(p) &
\left(\partial_{v_i}\partial_{w_j}\,|\hat\phi_\lambda|^2\right)(p) \\
\left(\partial_{w_i}\partial_{v_j}\,|\hat\phi_\lambda|^2\right)(p) &
\left(\partial_{w_i}\partial_{w_j}\,|\hat\phi_\lambda|^2\right)(p)
\endpmatrix \\
& =
|\hat\phi_\lambda|^2\, \pmatrix
\operatorname{Re}\left(\partial_{x_i}\partial_{x_j}\,\log\hat\phi_\lambda\right)(p) &
-\operatorname{Im}\left(\partial_{x_i}\partial_{x_j}\,\log\hat\phi_\lambda\right)(p) \\
-\operatorname{Im}\left(\partial_{x_i}\partial_{x_j}\,\log\hat\phi_\lambda\right)(p) &
-\operatorname{Re}\left(\partial_{x_i}\partial_{x_j}\,\log\hat\phi_\lambda\right)(p)
\endpmatrix\\
& =
|\hat\phi_\lambda|^2\, \pmatrix
\operatorname{Re}\operatorname{Hess}(\log\hat\phi_\lambda)(p) &
-\operatorname{Im}\operatorname{Hess}(\log\hat\phi_\lambda)(p) \\
-\operatorname{Im}\operatorname{Hess}(\log\hat\phi_\lambda)(p) &
-\operatorname{Re}\operatorname{Hess}(\log\hat\phi_\lambda)(p)
\endpmatrix.
\endalign
$$

The following fact is elementary.

\proclaim{Lemma 5.1} Assume that all critical points of $\hat\phi_\lambda$ in
$\hat Y$ are
non--degenerate. Then all critical points of $|\hat\phi_{\lambda}|^2$ in $\hat
Y$ are also
non--degenerate and have index equal to $n$.
\endproclaim

\demo{Proof}  One first observes that the characteristic values of the bilinear
form
$\operatorname{Hess}(|\hat\phi_{\lambda}|^2)(p)$ at a critical point $p$ occur
in pairs of
opposite sign. This follows from the immediate fact that, if
$\left(\smallmatrix a\\b\endsmallmatrix\right)$ is a characteristic vector with
characteristic value
$\alpha$, then $\left(\smallmatrix -b\\a\endsmallmatrix\right)$ is a
characteristic vector with
value $-\alpha$.  The result is thus proven if we show that the null space of
$\operatorname{Hess}(|\hat\phi_{\lambda}|^2)(p)$ is empty. We shall abbreviate
$\operatorname{Hess}(\log\hat\phi_{\lambda})(p)=H$,
$\operatorname{Re}H=R$ and $\operatorname{Im}H=I$. By contradiction, let us
suppose that 0 is a
characteristic value of $\operatorname{Hess}(|\hat\phi_{\lambda}|^2)(p)$. Since
by
assumption $\hat\phi_\lambda$ is nowhere zero on $\hat Y$, this means that
there exists a
non--zero real $2n$ vector $\left(\smallmatrix a\\b\endsmallmatrix\right)$ so
that $Ra=Ib$ and
$Ia=-Rb$. Thus, in particular, $Ra +\1 I a=Ib-\1 Rb$, i.e., $H(a+\1 b)=0$. But
by assumption $p$ is a non-degenerate critical point of $\hat\phi_{\lambda}$,
hence also of
$\log\hat\phi_{\lambda}$; so $H$ is non--singular. It follows that $a+\1 b=0$,
and---$a$,
$b$ being real---that $a=b=0$, a contradiction.
\enddemo

In view of Proposition 1.1, over the open dense subset $V-A\cap V$ of $\Lambda$
the critical points of $\phi_\lambda$ are all non--degenerate and constant in
number.
We may henceforth choose $\lambda$ to be a point in $V_{\Bbb Z^*}=(V-A\cap
V)\cap\Bbb
Z^{N-1}$, the subset of $V$ where all $\hat\lambda_i$ have vanishing imaginary
part and
non--zero integral real part. The hypersurface $\hat D$ is thus the support of
the divisor of the
global meromorphic function $\hat\phi_{\lambda}$, $\hat D=\hat D_0\cup\hat
D_\infty$, where
$\hat D_0=\cup_{\hat\lambda_i>0}\, \hat D_i$,
$\hat D_\infty=\cup_{\hat\lambda_i<0}\,\hat D_i$ are, respectively, the zero
and the polar locus
of $\hat\phi_\lambda$. Obviously, neither $\hat D_0$ nor $\hat D_\infty$ may be
empty. If one
were to directly apply standard Morse theory in the present set up, however,
one would encounter an
obstacle in the existence of points of indeterminacy for the function
$|\hat\phi_\lambda|^2$, that
is the points in
$\hat D_0\cap\hat D_\infty$, where $|\hat\phi_\lambda|^2$ has no limit. This
difficulty
admits a standard resolution which consists of recursively further blowing--up
$\hat X$ along the
components of the indeterminacy locus. There in fact exists (see e.g., \cite{4:
Section 2 of
Chapter 4}) a blow--up $X'@>\sigma'>>\hat X$ of $\hat X$ such that:
\roster
\item The supports of, respectively, the divisors of zeroes and of poles,
$D_0^\prime$ and
$D_\infty^\prime$, of the pull--back ${\sigma'}^*\hat\phi_\lambda$, are
disjoint;
\item $D'={\sigma'}^{-1}(\hat D)=D_0^\prime\cup D_\infty^\prime\cup D''$, the
function ${\sigma'}^*\hat\phi_\lambda$ being defined and nowhere vanishing on
the components of
$D''- D''\cap\bigl(D_0^\prime\cup D_\infty^\prime\bigr)$.
\endroster

Next, we shall ascertain that the pull-back ${\sigma'}^*\hat\phi_\lambda$, and
hence
${\sigma'}^*|\hat\phi_\lambda|^2$, has no critical points on
$D''- D''\cap\bigl(D_0^\prime\cup D_\infty^\prime\bigr)$.

\proclaim{Lemma 5.2} For $\lambda\in V_{\Bbb Z^*}$,
${\sigma'}^*|\hat\phi_\lambda|^2$ extends to a
positive $C^\infty$ function $F$ on $X'- D_\infty^\prime$ vanishing precisely
along $D_0^\prime$ and approaching infinity near $D_\infty^\prime$. Moreover,
this function has no critical points on
$D''- D''\cap\bigl(D_0^\prime\cup D_\infty^\prime\bigr)$.
\endproclaim

\demo{Proof} The first part of the statement---which we have included for
completeness---is true by
definition. We only have to examine the assertion about the critical points of
$F$.
Let $U$ be a small neighborhood in $\hat X$ of the intersection points of
two components of $\hat D$ along which the orders of $\hat\phi_\lambda$ are
opposite in
sign. Then, after a finite sequence $\pi_1$ of blow--ups along the
indeterminacy loci, one
arrives at the local situation where the only indeterminacy points of
$\pi_1^*\hat\phi_\lambda$ on
$\pi_1^{-1}(U)$ are those lying on the intersection of two irreducible
divisors $C_1=\{z_1=0\}$, $C_2=\{z_2=0\}$, the local form of
$\pi_1^*\hat\phi_\lambda$ being $z_1^{-m}\,z_2^m\,h(z_1,z_2,\dots)$, with some
positive non--zero
integer
$m$ and a nowhere vanishing holomorphic function $h$. After one more blow--up
$\pi_2$ along
$C_1\cap C_2$, given by $z_1=t_1,\,z_2=t_1t_2$ away from $C_1$, and denoting by
$\sigma'$ the composite
$\pi_2\pi_1$, one sees from the local form $t_2^m\,h(t_1,t_1t_2,\dots)$ that
the pull--back
${\sigma'}^*\hat\phi_\lambda$ is indeed holomorphic and non--zero along the
points on the
exceptional divisor $E=\{t_1=0\}$ which do not intersect the proper transform
of $C_1\cup
C_2$.  Let us now consider the derivatives of ${\sigma'}^*\hat\phi_\lambda$. In
particular, locally
near $E$ one has
$\partial_{t_2}\,{\sigma'}^*\hat\phi_\lambda=mt_2^{m-1}\,h(t_1,t_1t_2,\dots)+
t_2^m\,t_1\,\partial_{z_2}\,h(t_1,z_2,\dots)$, and
$\left.\partial_{t_2}\,{\sigma'}^*\hat\phi_\lambda\right|_E=m\,t_2^{m-1}\,h(0,0,\dots)$, which may
vanish only if $t_2=0$. It follows that ${\sigma'}^*\hat\phi_\lambda$, hence
also
${\sigma'}^*|\hat\phi_\lambda|^2$, has no critical points on the complement in
$E$ of the proper
transform of $C_2$. In other words, the only critical points, if any, of $F$ in
$D'$
necessarily lie on $D_0^\prime$, and the lemma has been proven.
\enddemo

Below is our main lemma. We shall henceforth denote by $\gamma$ the number of
critical points of $\hat\phi_\lambda$ in $\hat Y$ ($\lambda\in V-A\cap V$).
Note that, by Lemma
5.2, this is equal to the number of critical points of $F$ in $X'-
D_0^\prime\cup D_\infty^\prime$.

\proclaim{Theorem 5.3} Let $\lambda\in V_{\Bbb Z^*}$. Let also $\partial
\Ov{T}(D_0^\prime)$ be the boundary of an infinitesimally small (closed)
neighborhood
$\Ov{T}(D_0^\prime)$ of $D_0^\prime$ in $X'$. Then $X'- D_0^\prime\cup
D_\infty^\prime$ has the
homotopy type of $\partial\Ov{T}(D_0^\prime)$ with a number $\gamma$ of
$n$--cells attached.
\endproclaim

\demo{Proof} For $0<a<b$, let $\Phi$ denote the restriction of $F$ to
$M=F^{-1}[a,b]$. Since $D_0^\prime\cap\hat D_\infty^\prime$ is empty, if
neither $a$ or $b$ are critical values, $M$ is a compact real submanifold of
$X'-
D_0^\prime\cup D_\infty^\prime$ with smooth and disjoint boundary components
$\Phi^{-1}(a)$ and
$\Phi^{-1}(b)$. By Lemma 5.1 and Lemma 5.2,
$\Phi$ is a Morse function all of whose critical points have index $n$. Let
$\gamma(a,b)$ be the
number of critical points of $\Phi$ whose critical values lie in the interval
$(a,b)$.  In view of
the basic result of Morse theory (see for example \cite{5: Theorem 3.1 of
Chapter 6}), $M$ has the
homotopy type of  $\Phi^{-1}(a)$ with a number $\gamma(a,b)$ of $n$--cells
attached. If one chooses
$a$ and $b$ to be respectively so small and so large that $(a,b)$ contains all
critical values of
${\sigma'}^*|\hat\phi_\lambda|^2$, then, symbolically,
$$
M\,{\overset\text{hom}\to\cong}\,\Phi^{-1}(a)\cup e_1\cup\dots\cup e_\gamma,
$$
where the $e_i$ are $n$--cells. But clearly $\Phi^{-1}(a)$ is homotopic to
$\partial\Ov{T}(D_0^\prime)$;  on the other hand $M$ is homotopic to
$X'-T(D_0^\prime)\cup
T(D_\infty^\prime)$---where $T(D_\infty^\prime)$ is a small open neighborhood
of
$D_\infty^\prime$---and hence to $X'- D_0^\prime\cup D_\infty^\prime$.
\enddemo

\proclaim{Corollary 5.4} One has $\gamma=(-1)^n\,\chi(X'- D_0^\prime\cup
D_\infty^\prime)$.
\endproclaim

\demo{Proof} On the level of Euler characteristics, Theorem 5.3 implies the
relation
$$
\chi(X'- D_0^\prime\cup
D_\infty^\prime)=\chi\bigl(\partial\Ov{T}(D_0^\prime)\bigr)+(-1)^n\,\gamma.
$$
But $\partial\Ov{T}(D_0^\prime)$ is homotopic to the deleted neighborhood
$T(D_0^\prime)-D_0^\prime$, whose Euler characteristic is vanishing.
\enddemo

In order to complete the proof of Theorem 1.2 there only remains to observe
that
$$
\chi(Y)=\chi(X'- D_0^\prime\cup D_\infty^\prime).
$$
Since $Y$ is isomorphic to $\hat Y$ and $\hat Y$ to ${\sigma'}^{-1}(\hat Y)=X'-
D'$, the sought for
equality, $\chi(X'-D')=
\chi(X'- D_0^\prime\cup D_\infty^\prime)$, is equivalent---by the additivity
of the Euler characteristic---to the fact that
$\chi\bigl(D''- D''\cap (D_0^\prime\cup D_\infty^\prime)\bigr)=0$. One deduces
from the
explicit description of $\sigma'$ given in the proof of Lemma 5.2 that $D''$ is
a disjoint
union of exceptional divisors. With the same notation used above, the component
$E$ of $D''$
is by definition the projectivization of the normal bundle of $C_1\cap C_2$.
Let $\tilde C_1$, $\tilde C_2$ be respectively the proper transforms of $C_1$,
$C_2$ in $X'$. The complement $E-E\cap(\tilde C_1\cap\tilde C_2)$ is thus a
fiber bundle over
$\tilde C_1\cap\tilde C_2$ with fiber isomorphic to $\Bbb C^*=\Bbb C-\{0\}$. It
follows that
$\chi\bigl(E-E\cap(\tilde C_1\cap\tilde C_2)\bigr)=0$, and, summing over the
various components,
also that $\chi\bigl(D''- D''\cap (D_0^\prime\cup D_\infty^\prime)\bigr)=0$, as
desired. This concludes our second proof of Theorem 1.2.

\definition{Example 5.5} The operation of attaching a cell $e$ with boundary
$\dot e$ to a
topological space $S$ consists of providing an attaching map $s\:\dot e\to S$
and of
identifying every $x\in\dot e$ with $s(x)$.
The content of Theorem 5.3 is illustrated by the following simplest example of
Theorem 1.2.
Let $D=\{t_1,\dots,t_{N-1},t_N=\infty\}$ be a set of distinct points in $\PC$,
and let
$\phi_\lambda(x)=\prod_{i=1}^{N-1}\,(x-t_I)^{\lambda_I}$. Here
$N\ge 2$ and $\lambda=(\lambda_1,\dots,\lambda_N)$ is a point on the hyperplane
$\Lambda=\{\lambda_1+\dots+\lambda_N=0\}\subset\Bbb C^N$. One easily verifies
that, for a generic
$\lambda$ in $\Lambda$, all critical points of $\phi_\lambda$ are
non--degenerate. Moreover, if
$A=\{\lambda_1\cdots\lambda_N=0\}$ is the union of the coordinate hyperplanes
in $\Bbb C^N$, the
number of critical points for a generic
$\lambda\in\Lambda-A\cap\Lambda$ is equal to $N-2$. One may choose $\lambda$ so
that
$D_\infty=\{\infty\}$, $D_0=\{t_1,\dots,t_{N-1}\}$. Theorem 4.3 says in this
case that $\Bbb C$ is,
homotopically, the space obtained by attaching $N-2$ open segments
$e_i=(t_i,t_{i+1})$ to the
points in $D_0$, which we may assume to be ordered as $\operatorname{Re}t_i\le
\operatorname{Re}t_{i+1}$.

\enddefinition

\enddocument